# Uncertainty Marginal Price, Transmission Reserve, and Day-ahead Market Clearing with Robust Unit Commitment

Hongxing Ye, *Member, IEEE,* Yinyin Ge, Mohammad Shahidehpour, *Fellow, IEEE,* Zuyi Li, *Senior Member, IEEE*

*Abstract*—The increasing penetration of renewable energy in recent years has led to more uncertainties in power systems. These uncertainties have to be accommodated by flexible resources (i.e. upward and downward generation reserves). In this paper, a novel concept, Uncertainty Marginal Price (UMP), is proposed to price both the uncertainty and reserve. At the same time, the energy is priced at Locational Marginal Price (LMP). A novel market clearing mechanism is proposed to credit the generation and reserve and to charge the load and uncertainty within the Robust Unit Commitment (RUC) in the Day-ahead market. We derive the UMPs and LMPs in the robust optimization framework. UMP helps allocate the cost of generation reserves to uncertainty sources. We prove that the proposed market clearing mechanism leads to partial market equilibrium. We find that transmission reserves must be kept explicitly in addition to generation reserves for uncertainty accommodation. We prove that transmission reserves for ramping delivery may lead to Financial Transmission Right (FTR) underfunding in existing markets. The FTR underfunding can be covered by congestion fund collected from uncertainty payment in the proposed market clearing mechanism. Simulations on a six-bus system and the IEEE 118-bus system are performed to illustrate the new concepts and the market clearing mechanism.

*Index Terms*—Uncertainty Marginal Price, Cost Causation, Robust Unit Commitment, Financial Transmission Right, Generation Reserve, Transmission Reserve

## Nomenclature

**Indices**

| | |
|---|---|
| $i, l, t$ | indices for generators, lines, and time intervals |
| $m, n$ | index for buses |
| $m_i$ | index of bus where unit $i$ is located |
| $k$ | index of the worst point for uncertainty |

**Functions and sets**

| | |
|---|---|
| ˆ | symbol for the optimal value of a variable |
| $\mathcal{F}$ | feasible set for UC and dispatch |
| $\mathcal{U}$ | uncertainty set |
| $C_i^P(\cdot), C_i^I(\cdot)$ | cost related to dispatch and UC for unit $i$ |
| $\mathcal{L}(\cdot)$ | Lagrangian function |
| $\mathcal{G}(m)$ | set of units located at bus $m$ |
| $\mathcal{K}$ | set of the indices for $\hat{\epsilon}^k$ |

$\mathcal{K}_{m,t}^{\text{up}}, \mathcal{K}_{m,t}^{\text{down}}$    set of indices $k$ for upward and downward UMPs at bus $m$ time $t$

**Constants**

| | |
|---|---|
| $N_D, N_T$ | number of buses and time intervals |
| $d_{m,t}$ | aggregated equivalent load |
| $F_l$ | transmission line flow limit |
| $\Gamma_{l,m}$ | shift factor for line $l$ with respect to bus $m$ |
| $P_i^{\min}, P_i^{\max}$ | minimum and maximum generation outputs |
| $r_i^u, r_i^d$ | ramping-up/down limits between sequential intervals |
| $R_i^u, R_i^d$ | ramping-up/down limits for uncertainty accommodation |
| $u_{m,t}$ | bound for uncertainty |
| $\hat{\epsilon}$ | $\hat{\epsilon}^k$ is the $k^{th}$ worst uncertainty vector in $\mathcal{K}$, $\hat{\epsilon}^k \in \mathbb{R}^{N_D N_T}$, $\hat{\epsilon}_{m,t}^k \in \mathbb{R}$ |
| $\text{FTR}_{m \rightarrow n}$ | FTR amount from bus $m$ to $n$ |

**Variables**

| | |
|---|---|
| $I_{i,t}$ | unit on/off status indicators |
| $y_{i,t}, z_{i,t}$ | unit start-up and shut-down indicators |
| $P_{i,t}$ | generation dispatch |
| $P_{m,t}^{\text{inj}}$ | net power injection |
| $\epsilon_{m,t}$ | uncertainty at bus $m$ time $t$ |
| $\mathcal{Z}$ | optimal value of problem (SP) given $(x, y, z, I, P)$ |
| $\Delta P_{i,t}$ | generation re-dispatch |
| $\Delta f_{l,t}^{\text{pos}}$ | transmission capacity reserve in positive direction |
| $\Delta f_{l,t}^{\text{neg}}$ | transmission capacity reserve in negative direction |
| $\Delta P_{m,t}^{\text{inj}}$ | net power injection change |
| $\lambda_t, \lambda_t^k$ | Lagrangian multipliers |
| $\alpha, \beta, \eta$ | non-negative Lagrangian multipliers |
| $\pi_{m,t}$ | marginal prices. $\pi_{m,t}^e$ for energy price; $\pi_{m,t}^{u,k}$ is the UMP for $k$th uncertainty point; $\pi_{m,t}^{u,\text{up}}$ for upward UMP; $\pi_{m,t}^{u,\text{down}}$ for downward UMP |
| $Q_{i,t}^{\text{up}}, Q_{i,t}^{\text{down}}$ | upward and downward generation reserves |
| $\Psi_{m,t}$ | charge for uncertainty source |
| $\Theta_{i,t}^G, \Theta_{l,t}^T$ | credits to generation reserve for unit $i$ and transmission reserve for line $l$ at time $t$ |

This work is supported by the U.S. National Science Foundation Grant ECCS-1549937. The early version of this work was available on arXiv July 06, 2015, titled "Market Clearing for Uncertainty, Generation Reserve, and Transmission Reserve". The authors are with the Galvin Center for Electricity Innovation at Illinois Institute of Technology, Chicago, IL 60616, USA. (e-mail: hye@hawk.iit.edu; yge9@hawk.iit.edu; ms@iit.edu; lizu@iit.edu).

## I. Introduction

IN modern power systems, uncertainties grow significantly with the increasing penetration of Renewable Energy Source (RES), such as wind power generation. They pose new challenges for the operation of electricity markets. In the Day-ahead market (DAM), the Unit Commitment (UC) and Economic Dispatch (ED) problems considering uncertainties have become a focus of research in recent years. The objective of the UC problem is to find the least cost UC solution for the second day while respecting both system-wide and unit-wise constraints. By fixing the UC variables, the ED problem is established. The Locational Marginal Price (LMP) and reserve price are then obtained as byproducts of the ED problem [1], [2]. When considering the uncertainties, the generation from uncontrollable RES are uncertain parameters in the optimization problem.

Recently, Robust UC (RUC) is proposed to address the issues of uncertainty [3]–[7]. The largest merit is that the UC solution can be immunized against all the uncertainties in predefined set. The key idea of the two-stage RUC is to determine the optimal UC in the first stage which leads to the least cost for the worst scenario in the second stage. However, this approach is conservative and the Robust ED (RED) is absent. Authors in [8] combined the stochastic and robust approach using a weight factor in the objective function to address the conservativeness issue. [9], [10] employed the Affine Policy (AP) to formulate and solve the RED problem. A Multi-stage RUC is proposed to incorporate the latest information in each stage [11], where AP is also used to overcome the computational challenge. Recently, we reported a new approach which tries to bridge the gap of RUC and RED [7], [12].

In DAM, the main difficulty for pricing is that RED is absent in the traditional RUC [4], [5], [8]. On the other hand, a large number of works on pricing reserves exists within the UC framework considering contingencies [2], [13], [14] and stochastic security [15]. They are normally modeled as co-optimization problem. In [2], the reserve is cleared on zonal levels. Instead of countable contingency scenarios [15] or single additional scenario for reserve [16], the infinite continuous uncertainties are considered in the RUC, and the reserves are fully deliverable in infinite scenarios. In this paper, we propose a novel mechanism to price the energy, uncertainty, and flexibility simultaneously based on the RUC in [7]. An explicit price signal is derived for pricing the uncertainty. As the ED solution obtained is robust [7], both marginal impacts of the uncertainty and flexibility are reflected in these prices. In the proposed mechanism, reserve costs are allocated to uncertainty sources. Generation reserves, also called flexibilities in this paper, are the key factor for the robust optimization approaches. They are entitled to proper credits based on their contribution to uncertainty management. According to the market equilibrium analysis, market participants (price takers) can get the maximal profit by following the ISO/RTO's dispatch instruction.

The generation reserve and its deliverability are the main focus in [7]. The definition of LMP in [17] are employed to derive the energy price. The new concept, Uncertainty Marginal Price (UMP), is proposed to define the marginal cost of immunizing the next increment of uncertainty at a specific location. Load and generation are a pair, and they are priced at LMP. Uncertainty and flexibility (i.e., generation reserve) are another pair, and they are priced at UMP. Both LMPs and UMPs may vary with the locations due to transmission congestions. Limited by the transmission capacity and power flow equations, sometimes the uncertainties at certain buses cannot be mitigated by the system-wide cheapest generation reserve, and expensive generation reserve, which is deliverable, has to be kept in the system. Therefore, uncertainty sources are charged and generation reserves are credited based on UMPs at the corresponding buses.

As the transmission reserve is kept within the RUC framework, the congestion component may exist in both the energy price and reserve price even if the physical limit of the line is not reached yet in the base case scenario. LMP congestion costs are distributed to Financial Transmission Right (FTR) holders in the existing market according to the LMP difference and the FTR amount. The revenue inadequacy occurs when the LMP congestion cost collected is smaller than the credit distributed to FTR holders, which is also called FTR underfunding. This has been a serious issue in recent years in the industry [18], [19]. We reveal that transmission reserve will be another reason for FTR underfunding when physical transmission limit is adopted in Simultaneous Feasibility Test (SFT) for FTR market [18], [20], [21]. This conclusion is applicable to any robust UC framework for DAM.

The main contributions of this paper are listed as follows.

1) The novel UMP for uncertainties and generation reserves, as well as LMP for energy, are derived within a robust UC framework. The derivation is for uncertainties set with interval and budget constraints. The general concepts still apply when other uncertainty sets are modeled.
2) It is revealed that transmission capacities have to be reserved for uncertainty accommodation and the transmission reserves may cause FTR underfunding because of the deficiency of energy congestion revenues based on existing market rules.
3) A new market clearing mechanism is proposed to credit the generation and reserve and to charge the load and uncertainty. The payment collected from uncertainty sources can exactly cover the credits to generation reserves and transmission reserves, effectively resolving the FTR underfunding issue.

The rest of this paper is organized as follows. Derivation of the LMP and UMP is presented in Section II, so is the market clearing mechanism for charge and credit based on LMP and UMP. Case studies are presented in Section III. Section IV concludes this paper.

## II. RUC AND MARKET CLEARING

One motivation of this work is to price the uncertainty, and allocate the cost of uncertainty accommodation to the uncertainty source. As the uncertainty source is charged the uncertainty payment, it has the incentive to reduce the uncertainty. With UMP, we can follow the **cost causation principle**, which is normally required in the market design, to charge the uncertainty sources. Cost causation principle is described as "require that all approved rates reflect to some degree the costs actually caused by the customer who must pay them" in KN Energy, Inc. V. FERC, 968 F.2d 1295, 1300 (D.C. Cir. 1992).

Another important motivation is to provide a theory that supports the application of the RUC in the DAM clearing. Although the RUC/RED are studied extensively, the only application of the RUC now is for the Reliability Assessment Commitment (RAC) in the DAM. There are several reasons why they are not applied in the DAM clearing. First, the computation burden of RUC is much larger than the standard UC. Second, as the objective is the cost of the worst-case scenario [3], [5], the solution is criticized on over conservatism. Third, no economic dispatch and prices are available within the RUC framework. Recently, with the new achievements in the algorithms, models, and high-performance computing application [6]–[8], [22]–[24], the first two obstacles are being addressed with great promises. This paper tries to clear the last obstacle with the new model [7]. Adopting RUC in the market clearing can give clear price signals for the uncertainties and reserves. On the other side, it is also easier for the solution to pass the robustness test, which is a RUC, in RAC. To our best knowledge, this is the first work on pricing energy, uncertainties, and reserves within the robust optimization framework in DAM. Hence, we focus on illustrating the concept with the following assumptions.

- Network loss is ignored. Shift factor matrix is constant.
- Uncertainty is from load and RES. Contingency is ignored.
- The uncertainty budget set can be truly formulated by the ISO/RTO.

### A. RUC and RED

ISOs/RTOs desire to get the optimal UC and ED solution in the base-case scenario. They can re-dispatch the flexible resources, such as adjustable load demands and generators with fast ramping capabilities, to follow the load when deviation occurs (or uncertainty is revealed). Consistent with the robust literature [4], [5], the uncertainty set is modeled as

$$\mathcal{U} := \{ \boldsymbol{\epsilon} \in \mathbb{R}^{N_D N_T} : -u_{m,t} \leq \epsilon_{m,t} \leq u_{m,t}, \forall m, t$$
$$\sum_m \frac{|\epsilon_{m,t}|}{u_{m,t}} \leq \Lambda_t^\Delta, \forall t \}$$

$\Lambda_t^\Delta$ is the budget parameter and assumed as an integer [3]. It is noted that all the flexible resources are modeled as generators. In this paper, the RUC is formulated according to the model in [7].

$$\text{(RUC)} \min_{(\boldsymbol{x},\boldsymbol{p}) \in \mathcal{F}} C^I(\boldsymbol{x}) + C^P(\boldsymbol{p})$$

$$\text{s.t. } \boldsymbol{Ax} + \boldsymbol{Bp} \leq \boldsymbol{b} \tag{1}$$

$$\mathcal{F} := \Big\{ (\boldsymbol{x},\boldsymbol{p}) : \forall \boldsymbol{\epsilon} \in \mathcal{U}, \exists \Delta \boldsymbol{p} \text{ such that}$$

$$\boldsymbol{Cx} + \boldsymbol{Dp} + \boldsymbol{G}\Delta \boldsymbol{p} \leq \boldsymbol{d} + \boldsymbol{E}\boldsymbol{\epsilon} \Big\}. \tag{2}$$

The basic idea of the above model is to find a robust UC and ED for the base-case scenario. The UC $\boldsymbol{x}$ and dispatch $\boldsymbol{p}$ are immunized against any uncertainty $\boldsymbol{\epsilon} \in \mathcal{U}$. When uncertainty $\boldsymbol{\epsilon}$ occurs, it is accommodated by the generation adjustment $\Delta \boldsymbol{p}$. Please refer to Appendix A for the detailed formulation.

Column and Constraint Generation (CCG) based method is used to solve the above model [6]. Problem (MP) and (SP) are established as follows.

$$\text{(MP)} \quad \min_{(\boldsymbol{x},\boldsymbol{p})} C^I(\boldsymbol{x}) + C^P(\boldsymbol{p})$$

$$\text{s.t.} \quad \boldsymbol{Ax} + \boldsymbol{Bp} \leq \boldsymbol{b}$$
$$\boldsymbol{Cx} + \boldsymbol{Dp} + \boldsymbol{G}\Delta \boldsymbol{p}^k \leq \boldsymbol{d} + \boldsymbol{E}\hat{\boldsymbol{\epsilon}}^k, \forall k \in \mathcal{K} \tag{3a}$$

and

$$\text{(SP)} \quad \mathcal{Z} := \max_{\boldsymbol{\epsilon} \in \mathcal{U}} \min_{(\boldsymbol{s}, \Delta \boldsymbol{p}) \in \mathcal{R}(\boldsymbol{\epsilon})} \mathbf{1}^\top \boldsymbol{s} \tag{4a}$$

$$\mathcal{R}(\boldsymbol{\epsilon}) := \Big\{ (\boldsymbol{s}, \Delta \boldsymbol{p}) : \boldsymbol{s} \geq \mathbf{0} \tag{4b}$$

$$\boldsymbol{G}\Delta \boldsymbol{p} - \boldsymbol{s} \leq \boldsymbol{d} - \boldsymbol{Cx} - \boldsymbol{Dp} + \boldsymbol{E}\boldsymbol{\epsilon} \Big\} \tag{4c}$$

where $\mathcal{K}$ is the index set for uncertainty points $\hat{\boldsymbol{\epsilon}}$ which are dynamically generated in (SP) with iterations. Please refer to Appendix B for the detailed formulation. It should be noted that $\hat{\boldsymbol{\epsilon}}^k$ is the extreme point of $\mathcal{U}$. Variable $\Delta \boldsymbol{p}^k$ is associated with $\hat{\boldsymbol{\epsilon}}^k$. The objective function in (SP) is to find the worst point in $\mathcal{U}$ given $(\boldsymbol{x}, \boldsymbol{p})$. The procedure is

1: $\mathcal{K} \leftarrow \emptyset, k \leftarrow 1, \mathcal{Z} \leftarrow +\infty$, define feasibility tolerance $\delta$
2: **while** $\mathcal{Z} \geq \delta$ **do**
3:     Solve (MP), obtain optimal $(\hat{\boldsymbol{x}}, \hat{\boldsymbol{p}})$.
4:     Solve (SP) with $\boldsymbol{x} = \hat{\boldsymbol{x}}, \boldsymbol{p} = \hat{\boldsymbol{p}}$, get solution $(\mathcal{Z}, \hat{\boldsymbol{\epsilon}}^k)$
5:     $\mathcal{K} \leftarrow \mathcal{K} \cup k, k \leftarrow k + 1$
6: **end while**

Once the procedure is converged, we also get the optimal UC and ED solution by solving (MP). Similar to traditional LMP calculation, we fix the binary variables as $\hat{\boldsymbol{x}}$. Then a convex linear programming problem (RED) can be formed as

$$\text{(RED)} \min_{P, \Delta P} \quad \sum_t \sum_i C_i^P(P_{i,t}) \tag{5}$$

s.t.

$$(\lambda_t) \quad \sum_i P_{i,t} = \sum_m d_{m,t}, \forall t, \tag{6a}$$

$$(\bar{\beta}_{i,t}) \quad P_{i,t} \leq \hat{I}_{i,t} P_i^{\max}, \forall i, t \tag{6b}$$

$$(\underline{\beta}_{i,t}) \quad -P_{i,t} \leq -\hat{I}_{i,t} P_i^{\min}, \forall i, t \tag{6c}$$

$$(\bar{\alpha}_{i,t}) \quad P_{i,t} - P_{i,t-1} \leq r_i^u (1 - \hat{y}_{i,t}) + P_i^{\min} \hat{y}_{i,t}, \forall i, t \tag{6d}$$

$$(\underline{\alpha}_{i,t}) \quad -P_{i,t} + P_{i,t-1} \leq r_i^d (1 - \hat{z}_{i,t}) + P_i^{\min} \hat{z}_{i,t}, \forall i, t \tag{6e}$$

$$(\bar{\eta}_{l,t}) \quad \sum_m \Gamma_{l,m} P_{m,t}^{\text{inj}} \leq F_l, \forall l, t \tag{6f}$$

$$(\underline{\eta}_{l,t}) \quad -\sum_m \Gamma_{l,m} P_{m,t}^{\text{inj}} \leq F_l, \forall l, t \tag{6g}$$

$$(\lambda_t^k) \quad \sum_i \Delta P_{i,t}^k = \sum_m \hat{\epsilon}_{m,t}^k, \forall t, \forall k \in \mathcal{K} \tag{7a}$$

$$(\bar{\beta}_{i,t}^k) \quad P_{i,t} + \Delta P_{i,t}^k \leq \hat{I}_{i,t} P_i^{\max}, \forall i, t, \forall k \in \mathcal{K} \tag{7b}$$

$$(\underline{\beta}_{i,t}^k) \quad -P_{i,t} - \Delta P_{i,t}^k \leq -\hat{I}_{i,t} P_i^{\min}, \forall i, t, \forall k \in \mathcal{K} \tag{7c}$$

$$(\bar{\alpha}_{i,t}^k) \quad \Delta P_{i,t}^k \leq R_i^u (1 - \hat{y}_{i,t}), \forall i, t, \forall k \in \mathcal{K} \tag{7d}$$

$$(\underline{\alpha}_{i,t}^k) \quad -\Delta P_{i,t}^k \leq R_i^d (1 - \hat{z}_{i,t+1}), \forall i, t, \forall k \in \mathcal{K} \tag{7e}$$

$$(\bar{\eta}_{l,t}^k) \quad \sum_m \Gamma_{l,m} (P_{m,t}^{\text{inj}} + \Delta P_{m,t}^{\text{inj}}) \leq F_l, \forall l, t, \forall k \in \mathcal{K} \tag{7f}$$

$$(\underline{\eta}_{l,t}^k) \quad -\sum_m \Gamma_{l,m} (P_{m,t}^{\text{inj}} + \Delta P_{m,t}^{\text{inj},k}) \leq F_l, \forall l, t, \forall k \in \mathcal{K} \tag{7g}$$

where (6a)-(6g) are the constraints for the base ED, and (7a)-(7g) are constraints for different extreme points in $\mathcal{U}$. (7a) denotes the load balance after re-dispatch. The generation adjustments respects capacity limits (7b)(7c) and ramping limits (7d)(7e). Network constraints are denoted by (7f)(7g). The $P_{m,t}^{\text{inj}}$ and $\Delta P_{m,t}^{\text{inj},k}$ are defined as

$$P_{m,t}^{\text{inj}} := \sum_{i \in \mathcal{G}(m)} P_{i,t} - d_{m,t}, \forall m, t,$$

and

$$\Delta P_{m,t}^{\text{inj},k} := \sum_{i \in \mathcal{G}(m)} \Delta P_{i,t}^k - \hat{\epsilon}_{m,t}^k, \forall m, t, k,$$

respectively.

### B. Marginal Prices

In this section, marginal prices for the energy, uncertainty, and generation reserve are derived based on the Lagrangian function. Denote the Lagrangian function for (RED) as $\mathcal{L}(P, \Delta P, \lambda, \alpha, \beta, \eta)$, which is shown in Appendix C. According to the definition of marginal price [17], the LMP for energy at bus $m$ is

$$\pi_{m,t}^{\text{e}} = \frac{\partial \mathcal{L}(P, \Delta P, \lambda, \alpha, \beta, \eta)}{\partial d_{m,t}} \tag{8}$$
$$= \lambda_t - \sum_l \Gamma_{l,m} \left( \bar{\eta}_{l,t} - \underline{\eta}_{l,t} \right) - \sum_l \sum_{k \in \mathcal{K}} \Gamma_{l,m} \left( \bar{\eta}_{l,t}^k - \underline{\eta}_{l,t}^k \right)$$

It is observed that the impact of the uncertainty is also reflected in the LMP.

The new concept, UMP for DAM, is defined as the marginal cost of immunizing the next unit increment of uncertainty. For $\hat{\epsilon}^k$, an extreme point of $\mathcal{U}$, the UMP is

$$\pi_{m,t}^{\text{u},k} = \frac{\partial \mathcal{L}(P, \Delta P, \lambda, \alpha, \beta, \eta)}{\partial \hat{\epsilon}_{m,t}^k} = \lambda_t^k - \sum_l \Gamma_{l,m} \left( \bar{\eta}_{l,t}^k - \underline{\eta}_{l,t}^k \right) \tag{9}$$

Both the uncertainty and generation reserve are priced at $\pi_{m,t}^{\text{u},k}$. In the derivation of $\pi_{m,t}^{\text{u},k}$, the worst point $\hat{\epsilon}^k$ is the only concern. Therefore, the general principles in this paper still work when $\mathcal{U}$ is replaced with other sets.

It should be noted that (9) is intermediate price signals. In order to get the aggregated UMPs, the following new sets are defined based on the sign of $\pi_{m,t}^{\text{u},k}$.

$$\mathcal{K}_{m,t}^{\text{up}} := \{k : \pi_{m,t}^{\text{u},k} \geq 0\}; \quad \mathcal{K}_{m,t}^{\text{down}} := \{k : \pi_{m,t}^{\text{u},k} < 0\} \tag{10}$$

The aggregated upward and downward UMPs are defined as

$$\pi_{m,t}^{\text{u,up}} := \sum_{k \in \mathcal{K}_{m,t}^{\text{up}}} \pi_{m,t}^{\text{u},k}; \quad \pi_{m,t}^{\text{u,down}} := \sum_{k \in \mathcal{K}_{m,t}^{\text{down}}} \pi_{m,t}^{\text{u},k} \tag{11}$$

respectively. In the following context, we will show how the aggregated UMPs are used.

### C. Market Clearing Mechanism

With LMP and UMP, the charges and credits for the market participants become clear and fair in the DAM. Energy clearing is straightforward. The basic principle related to uncertainty and flexibility is that those who cause uncertainties (uncertainty sources), such as RES, pay based on UMP and those who contribute to the management of uncertainties (uncertainty mitigators), such as generators or storage with ramping capabilities, get paid.

*1) Energy Payment and Credit:* LSEs pay based on the amount of the load and LMP. The energy payment from the LSE at Bus $m$ at $t$ is $\pi_{m,t}^{\text{e}} d_{m,t}$. It should be noted that RES is entitled to the credit due to the negative load modeled in RUC. Generator $i$, located at Bus $m_i$, is entitled to the credit $\pi_{m_i,t}^{\text{e}} P_{i,t}$ for energy production.

*2) Charge to Uncertainty Source:* The uncertainty source can be charged as

$$\Psi_{m,t} = \sum_{k \in \mathcal{K}} \pi_{m,t}^{\text{u},k} \hat{\epsilon}_{m,t}^k \tag{12}$$

The uncertainty source pays based on the marginal price and the worst point $\hat{\epsilon}^k$. The uncertainty source is charged only when $\pi_{m,t}^{\text{u},k}$ is non-zero, and it may have to pay more when the uncertainty becomes larger. The uncertainty point $\hat{\epsilon}_{m,t}^k$ can be upward (i.e. $\hat{\epsilon}_{m,t}^k \geq 0$) or downward (i.e. $\hat{\epsilon}_{m,t}^k \leq 0$). We have the following lemma regarding the relation between the signs of $\pi_{m,t}^{\text{u},k}$ and $\hat{\epsilon}_{m,t}^k$.

**Lemma 1.** If $\hat{\epsilon}_{m,t}^k > 0$, then $\pi_{m,t}^{\text{u},k} \geq 0$. If $\hat{\epsilon}_{m,t}^k < 0$, then $\pi_{m,t}^{\text{u},k} \leq 0$.

Please check Appendix D-A for the proof. When the budget set is adopted, the extreme point $\hat{\epsilon}_{m,t}^k \in \{-u_{m,t}, 0, u_{m,t}\}$ [4], [7], so the uncertainty charge in (12) can also be written as (13) according to Lemma 1 and (11).

$$\Psi_{m,t} = \pi_{m,t}^{\text{u,up}} u_{m,t} + \pi_{m,t}^{\text{u,down}}(-u_{m,t}) \tag{13}$$

Thus, upward and downward uncertainties are charged separately. It should be noted that we still need to use (12) when other uncertainty sets are used.

*3) Credit to Generation Reserve:* Only resources that can provide deliverable generation reserve are entitled to credits. If $i \in \mathcal{G}(m)$, then the credits can be formulated as

$$\Theta_{i,t}^G = \sum_{k \in \mathcal{K}} \pi_{m,t}^{\text{u},k} \Delta P_{i,t}^k. \tag{14}$$

In other words, generation reserve is paid the UMP at the bus where it is located. If $\pi_{m,t}^{\text{u},k} = 0$, then the associated credit is zero no matter what the value of $\Delta P_{i,t}^k$ is. Similar to the uncertainties, the generation reserves can be in either upward or downward direction. Denote the upward generation reserve as $Q_{i,t}^{\text{up}}$ and the downward generation reserve as $Q_{i,t}^{\text{down}}$

$$Q_{i,t}^{\text{up}} := \min \left\{ \hat{I}_{i,t} P_i^{\max} - P_{i,t}, \ R_i^u (1 - \hat{y}_{i,t}) \right\}, \tag{15}$$
$$Q_{i,t}^{\text{down}} := \max \left\{ \hat{I}_{i,t} P_i^{\min} - P_{i,t}, \ -R_i^d (1 - \hat{z}_{i,t+1}) \right\}. \tag{16}$$

We also have the following lemma regarding the relation between $Q_{i,t}^{\text{up}}$, $Q_{i,t}^{\text{down}}$ and $\Delta P_{i,t}^k$.

**Lemma 2.** If $i \in \mathcal{G}(m)$, then the optimal solution $\Delta P_{i,t}^k$ to problem (RED) is

$$\Delta P_{i,t}^k = \begin{cases} Q_{i,t}^{\text{up}}, & \text{if } \pi_{m,t}^{u,k} > 0 \\ Q_{i,t}^{\text{down}}, & \text{if } \pi_{m,t}^{u,k} < 0 \end{cases}$$

and

$$\pi_{m,t}^{u,k} = \bar{\beta}_{i,t}^k - \underline{\beta}_{i,t}^k + \bar{\alpha}_{i,t}^k - \underline{\alpha}_{i,t}^k, \tag{17}$$

Please check Appendix D-B for the proof. The credit to generation reserve $i$ located at bus $m$ (14) can be rewritten as (18) according to Lemma 2 and (11).

$$\Theta_{i,t}^G = \pi_{m,t}^{u,\text{up}} Q_{i,t}^{\text{up}} + \pi_{m,t}^{u,\text{down}} Q_{i,t}^{\text{down}} \tag{18}$$

(18) shows that the upward and downward generation reserves are credited separately. Flexible resources may receive credits for both the upward and downward generation reserves simultaneously. (18) always holds even if other uncertainty sets are modeled in RUC.

### D. Transmission Reserve and Revenue Adequacy

Some transmission capacities are reserved according to the solution to (RED). These transmission reserves are used to ensure the ramping deliverability when the uncertainty is revealed, as shown in (7f) and (7g). It is noted that they are determined automatically in (RED), and kept explicitly without explicit transmission reserve requirement constraints. Just like the "scheduled" generation reserve, the "scheduled" transmission reserves in positive direction and negative direction are

$$\Delta f_{l,t}^{\text{pos}} := F_l - \sum_m \Gamma_{l,m} P_{m,t}^{\text{inj}}, \tag{19}$$

$$\Delta f_{l,t}^{\text{neg}} := F_l + \sum_m \Gamma_{l,m} P_{m,t}^{\text{inj}}, \tag{20}$$

respectively. They are always non-negative.

An important issue related to the transmission reserve is the credit entitled to the Financial Transmission Right (FTR) holders. FTR is a financial instrument used to hedge congestion cost in the electricity market, where participants are charged or credited due to the transmission congestion [21], [25]. Within the robust framework, the effective transmission capacity for base-case scenario is different from the physical limit, which is used in the Simultaneous Feasibility Test (SFT) for FTR market [18], [20], [21]. In the existing market, the FTR credit is funded by the energy congestion cost, which is the net payment of energy. However, the energy congestion cost may not be sufficient to fund the FTR credit [18], [19]. We argue that the transmission reserve becomes a new reason for FTR underfunding in any framework to guarantee the ramping deliverability.

**Theorem 1.** If transmission reserve $\Delta f_{l,t}^{\text{pos}}$ and $\Delta f_{l,t}^{\text{neg}}$ are kept for line $l$ at time $t$ in DAM, then the maximum FTR underfunding associated with line $l$ at time $t$ is

$$\sum_{k \in \mathcal{K}} \left( \bar{\eta}_{l,t}^k \Delta f_{l,t}^{\text{pos}} + \underline{\eta}_{l,t}^k \Delta f_{l,t}^{\text{neg}} \right) \tag{21}$$

due to the deficiency of energy congestion cost.

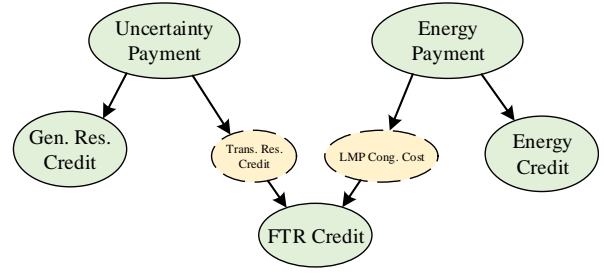

Fig. 1. Money flow of the proposed market clearing mechanism, where uncertainty sources make the uncertainty payment, and LSEs make the energy payment.

Please check Appendix D-C for the proof. From the FTR holder's point of view, (21) is the credit due to the transmission reserve. Therefore, we also call (21) transmission reserve credit, and denote it by

$$\Theta_{l,t}^T := \sum_{k \in \mathcal{K}} \left( \bar{\eta}_{l,t}^k \Delta f_{l,t}^{\text{pos}} + \underline{\eta}_{l,t}^k \Delta f_{l,t}^{\text{neg}} \right). \tag{22}$$

At most one of $\bar{\eta}_{l,t}^k$ and $\underline{\eta}_{l,t}^k$ is non-zero for transmission $l$. The credit to positive transmission reserve is zero for line $l$ at time $t$, when either $\sum_{k \in \mathcal{K}} \bar{\eta}_{l,t}^k = 0$ or $\Delta f_{l,t}^{\text{pos}} = 0$.

**Theorem 2.** If (RED) is feasible, then uncertainty payment can exactly cover generation reserve credit and transmission reserve credit, and the revenue adequacy is always guaranteed in the proposed market clearing mechanism.

Please check Appendix D-D for the proof. Theorem 1 reveals that FTR underfunding issue can occur within the existing market structures as long as the transmission reserve is non-zero, even if the LMPs are calculated based on other approaches. Theorem 2 shows that the new market clearing mechanism overcomes the FTR underfunding issue. The money flow of the proposed market clearing mechanism is depicted in Fig.1. Energy payment collected based on LMP is distributed to FTR holders as LMP congestion cost and generators as energy credit. On the other hand, the payment collected based on UMP is distributed to FTR holders as transmission reserve credit and flexible resources as generation reserve credit. The LMP congestion cost and transmission reserve credit can exactly cover the FTR credit, which is calculated based on the LMP difference and FTR amount.

### E. Market Equilibrium

In this section, we characterize the competitive market equilibrium model. In the electricity industry, the partial market equilibrium model is often employed [1], [13], [26], where market participants are price takers [27].

The energy is cleared according to (6a). Uncertainty and generation reserve are cleared according to (7a). Without loss of generality, consider unit $i$ located at bus $m$. Its profit maximization problem can be formulated as

$$(\text{PMP}_i) \max_{P_{i,t}} \sum_t \left( \begin{array}{c} P_{i,t} \pi_{m,t}^e + \pi_{m,t}^{u,\text{up}} Q_{i,t}^{\text{up}} + \pi_{m,t}^{u,\text{down}} Q_{i,t}^{\text{down}} \\ -C_i^P(P_{i,t}) \end{array} \right)$$

$$\text{s.t. } (6b) - (6e), (15) - (16)$$

where the decision variable is $P_{i,t}$ given the price signal $(\pi^{e}_{m,t}, \pi^{u,up}_{m,t}, \pi^{u,down}_{m,t})$. As proved in Appendix D-E, unit $i$ is not inclined to change its power output level as it can obtain the maximum profit by following the ISO's dispatch instruction $\hat{P}_{i,t}$. Price signal $\pi^{e}_{m,t}$ provides the incentives for unit $i$ to dispatch power output to $\hat{P}_{i,t}$, and price $\pi^{u,k}_{m,t}$ gives incentives for unit $i$ to maintain the generation reserve for uncertainty. Hence, the dispatch instruction $\hat{P}_{i,t}$ and price signal $(\pi^{e}_{m,t}, \pi^{u,up}_{m,t}, \pi^{u,down}_{m,t})$ constitute a competitive partial equilibrium [27].

*F. Discussions*

As the $P_{i,t}$ and $\Delta P^{k}_{i,t}$ (or $(Q^{up}_{i,t}, Q^{down}_{i,t})$) are coupled by (7b) and (7c), the opportunity cost $(\bar{\beta}^{k}_{i,t} - \underline{\beta}^{k}_{i,t})$ is enough to provide the incentives for $i$ to keep the generation level at $\hat{P}_{i,t}$. Including $\bar{\alpha}^{k}_{i,t} - \underline{\alpha}^{k}_{i,t}$ in the generation reserve price has several benefits. Firstly, generation reserves provided by different units are priced fairly. Generation reserve prices are the same for the units at the same bus, and they may vary with locations if line congestions exist. Secondly, higher generation reserve price attracts long-term investment for flexible resources. Thirdly, it is consistent with the existing reserve pricing practice [2], [28]. In fact, generation reserve price is consistent with the UMP. Therefore, the uncertainties and flexibilities are also treated fairly at the same bus.

The upward and downward UMPs are obtained according to (11), respectively. The uncertainty sources are charged according to (13). The generation reserves are credited according to (18). The price signal $\pi^{u,k}_{m,t}$ defined in (9) and re-dispatch $\Delta P^{k}_{i,t}$ are intermediate variables for market clearing. The proposed UMP may be non-zero even if the uncertainty at a bus is zero. This is similar to the LMP, which may also be non-zero for the bus without load.

The market clearing mechanism proposed in this paper follows the cost causation principle for the cost allocation. In reality, it may be controversial to allocate the reserve cost to uncertainty sources. However, we argue that it would be fair and must be done when the RES penetration level is high. An extreme case is when the loads are all supplied by RES. There has been study showing it is possible that the increasing RES penetration can cause higher system operation cost. This issue cannot be handled by the existing market clearing mechanism, in which loads pay for the additional system reserve that is required to accommodate the uncertainty from RES. In other words, loads are actually providing subsidies to RES. When the RES penetration level is low, the subsidies can help the growth of the RES. However, when the RES penetration level is high, these growing subsidies will cause serious fairness issue. On the other hand, with UMP as the stimulating price signals, RES will have the incentives to improve its forecast techniques and reduce its uncertainty. In the ideal case when its uncertainty approaches zero, RES will no longer pay.

Following the existing practice, the UC variables are fixed during the marginal price derivation. Hence, the uplift issue, which exists in the real market, still remains in the proposed market clearing mechanism. Although the UC variables are fixed, the LMP and reserve price in the real market can provide effective signals for the long-term investment of generation and transmission as well as consumption strategy of electricity. Similarly, the uncertainty impact is not only reflected in UC, but also in the ED within the RUC model in this paper. Hence, the proposed LMP and UMP can also provide signals for the long-term investment of flexibilities (i.e. generation, transmission, and demand).

The pricing for uncertainties proposed in this paper is not in conflict with the pricing for traditional reserves, which are mainly prepared for the contingencies. The traditional reserve prices can be derived in the framework by adding extra traditional reserve constraints, and the corresponding reserve costs can still be allocated to LSEs.

It is observed that the credit in (14) is the sum of credits for all extreme points. That is because the related constraints may be binding for multiple extreme points, and the dual variables (shadow prices) for these constraints work together in the dual problem. The traditional price for energy and reserve also has similar form when multiple contingencies are modeled.

Although only one scenario will happen in reality, we still need to consider the worst scenario defined in uncertainty set and keep enough reserves in DAM. That is because DAM is a financial market, and the LMP and UMP are the financially binding prices. This is similar to the existing market model considering contingencies. Even if the contingency seldom occurs, they are still modeled for market clearing, and the contingencies are reflected in LMP and reserve price.

The issue of price multiplicity still exists in the proposed model [29] because problem (RED) is a linear programming (LP) problem. However, the price is unique with the nonde-generacy assumption. For simplicity, we have considered a single-sided auction in the proposed model. By introducing the demand bids, we can formulate a double-sided auction and the general principles in this paper will still apply.

## III. CASE STUDY

A six-bus system and the IEEE 118-Bus system are simulated to illustrate the proposed market clearing mechanism. In the six-bus system, the basic ideas of UMP are presented within the robust optimization framework. FTR underfunding issue is illustrated and a comparison between the UMP and traditional reserve price are presented. In the IEEE 118-Bus system, the UMP related products are presented for different uncertainty levels. The behaviors and impacts of flexible sources are analyzed by an energy storage example.

### A. Six-bus System

A six-bus system is studied in this section. The one-line diagram is shown in Fig. 2. The unit data and line data are shown in Table I and Table II, respectively. Table III presents the load and uncertainty information. Column "Base Load" shows the hourly forecasted load. Assume that the load distributions are 20%, 40%, and 40% for Bus 3, Bus 4, and Bus 5, respectively. $\bar{u}_{1,t}$ and $\bar{u}_{3,t}$ in Table III are the bounds of the uncertainties at Bus 1 and Bus 3, respectively. The uncertainty bounds at other buses are 0.

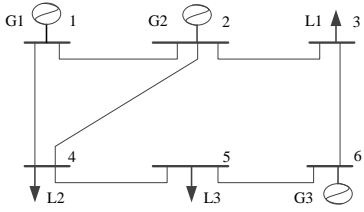

Fig. 2. One-line diagram for 6-bus system.

### TABLE I
#### Unit Data for the 6-bus System

| # | $P^{min}$ | $P^{max}$ | $P_0$ | a | b | c | $R^u$ | $R^d$ | $C_u$ | $C_d$ | $T^{on}$ | $T^{off}$ | $T_0$ |
|---|---|---|---|---|---|---|---|---|---|---|---|---|---|
| 1 | 100 | 220 | 120 | 0.004 | 13.5 | 176.9 | 24 | 24 | 180 | 50 | 4 | 4 | 4 |
| 2 | 10 | 100 | 50 | 0.001 | 32.6 | 129.9 | 12 | 12 | 360 | 40 | 3 | 2 | 3 |
| 6 | 10 | 20 | 0 | 0.005 | 17.6 | 137.4 | 5 | 5 | 60 | 0 | 1 | 1 | −2 |

$P^{min}, P^{max}, P_0$: min/max/initial generation level (MW);
fuel cost ($): $aP^2 + bP + c$ ;
$R^u, R^d$: ramping up/down rate (MW/h);
$C_u, C_d$: startup/shutdown cost ($);
$T^{on}, T^{off}, T_0$: min on/min off/initial time (h)

It is assumed that the relative forecasting errors increase with hours. Uncertainty $\epsilon_{1,t}$ and $\epsilon_{3,t}$ also respect

$$-\Lambda \cdot \bar{u}_{m,t} \le \epsilon_{m,t} \le \Lambda \cdot \bar{u}_{m,t}, \forall t, m \quad (23a)$$

$$\sum_m \frac{|\epsilon_{m,t}|}{\bar{u}_{m,t}} \le \Lambda_\Delta, \forall t, \quad (23b)$$

where (23a) denotes the uncertainty interval at a single bus, and (23b) represents the system-wide uncertainty [5]. The $\Lambda$ and $\Lambda^\Delta$ are the budget parameters for the single bus and system, respectively.

*1) LMP and UMP:* Consider the case where $\Lambda = 1, \Lambda^\Delta = 2$. The CCG based approach converges after 2 iterations.

### TABLE II
#### Line Data for the 6-bus System

| from | 1 | 1 | 2 | 5 | 3 | 2 | 4 |
|---|---|---|---|---|---|---|---|
| to | 2 | 4 | 4 | 6 | 6 | 3 | 5 |
| x(p.u.) | 0.17 | 0.258 | 0.197 | 0.14 | 0.018 | 0.037 | 0.037 |
| capacity(MW) | 200 | 100 | 100 | 100 | 100 | 200 | 200 |

### TABLE III
#### Load and Uncertainty Data for the 6-bus System (MW)

| Time (h) | Base Load | $\bar{u}_{1,t}$ | $\bar{u}_{3,t}$ | Time (h) | Base Load | $\bar{u}_{1,t}$ | $\bar{u}_{3,t}$ |
|---|---|---|---|---|---|---|---|
| 1 | 175.19 | 1.09 | 0.29 | 13 | 242.18 | 19.68 | 5.25 |
| 2 | 165.15 | 2.06 | 0.55 | 14 | 243.6 | 21.32 | 5.68 |
| 3 | 158.67 | 2.98 | 0.79 | 15 | 248.86 | 23.33 | 6.22 |
| 4 | 154.73 | 3.87 | 1.03 | 16 | 255.79 | 25.58 | 6.82 |
| 5 | 155.06 | 4.85 | 1.29 | 17 | 256 | 27.2 | 7.25 |
| 6 | 160.48 | 6.02 | 1.6 | 18 | 246.74 | 27.76 | 7.4 |
| 7 | 173.39 | 7.59 | 2.02 | 19 | 245.97 | 29.21 | 7.79 |
| 8 | 177.6 | 8.88 | 2.37 | 20 | 237.35 | 29.67 | 7.91 |
| 9 | 186.81 | 10.51 | 2.8 | 21 | 237.31 | 31.15 | 8.31 |
| 10 | 206.96 | 12.94 | 3.45 | 22 | 232.67 | 31.99 | 8.53 |
| 11 | 228.61 | 15.72 | 4.19 | 23 | 195.93 | 28.16 | 7.51 |
| 12 | 236.1 | 17.71 | 4.72 | 24 | 195.6 | 29.34 | 7.82 |

### TABLE IV
#### Marginal Costs at Different Generation Levels ($/MWh)

| Gen. 1 | | | Gen. 2 | | | Gen. 3 | | |
|---|---|---|---|---|---|---|---|---|
| $\underline{P}_1^w$ | $\bar{P}_1^w$ | mar. cost | $\underline{P}_2^w$ | $\bar{P}_2^w$ | mar. cost | $\underline{P}_3^w$ | $\bar{P}_3^w$ | mar. cost |
| 100 | 124 | 14.396 | 10 | 28 | 32.638 | 10 | 12 | 17.71 |
| 124 | 148 | 14.588 | 28 | 46 | 32.674 | 12 | 14 | 17.73 |
| 148 | 172 | 14.78 | 46 | 64 | 32.71 | 14 | 16 | 17.75 |
| 172 | 196 | 14.972 | 64 | 82 | 32.746 | 16 | 18 | 17.77 |
| 196 | 220 | 15.164 | 82 | 100 | 32.782 | 18 | 20 | 17.79 |

### TABLE V
#### Generation and Reserve ($\Lambda = 1, \Lambda^\Delta = 2$, MW)

| T | $P_1$ | $P_2$ | $P_3$ | $Q_1^{up}$ | $Q_1^{down}$ | $Q_2^{up}$ | $Q_2^{down}$ | $Q_3^{up}$ | $Q_3^{down}$ |
|---|---|---|---|---|---|---|---|---|---|
| 21 | 195.19 | 25.58 | 16.54 | 24 | −24 | 12 | −12 | 3.46 | −5 |

Hence, $\mathcal{K} = \{1, 2\}$. Given the UC solutions, the problem (RED) can be solved by commercial LP solver. The marginal prices are then obtained as byproducts.

The generation outputs are presented in Table V at Hours 21. It can be observed that G1 supplies most of the loads at Hour 21, which is 195.19 MW. According to the bid information in Table IV, G2 is much more expensive than G1 and G3. Hence, the output of G2 is relatively small and at the low level of its capacity. The upward and downward generation reserves provided by the three units are also listed in Table V. These data can be obtained directly from Eqs. (15) and (16) given the generation output $P_{i,t}$. Although the remaining generation capacity of G1 is $220 - 195.19 = 24.62$ MW, the upward reserve is limited by its upward ramping rate 24 MW. In the meantime, the upward reserve provided by G3 is limited by its generation capacity although it has more remaining ramping capacity (i.e. $\min\{20 - 16.54, 5\} = 3.46$ MW).

Table VI shows the extreme points obtained in the CCG-based approach. The intermediate price signals for these points $\pi_{m,t}^{u,k}$ are also presented. It can be observed that the worst point is always obtained at the extreme point of the uncertainty set. For example, at Hour 21, the $\hat{\epsilon}_{1,1}^1$ is 31.15 MW. It is exactly the upper bound of the uncertainty at Hour 21 at Bus 1. The data in Table VI also verifies Lemma 1. The intermediate UMPs $\pi_{m,t}^{u,k}$ have the same sign as the uncertainties $\hat{\epsilon}_{m,t}^k$ at the same bus.

The LMPs, aggregated upward UMPs, and aggregated downward UMPs at Hour 21 are shown in Table VII. It is noted that UMPs still exist at buses without uncertainties (i.e., Buses 2,4,5,6). This is similar to LMPs, which also exist at buses where net power injections are 0. The LMPs vary with locations, which indicates that the line congestion exists.

### TABLE VI
#### Extreme Points of Uncertainty Set

| | k = 1 | | | | k = 2 | | | |
|---|---|---|---|---|---|---|---|---|
| $t$ | $\hat{\epsilon}_{1,t}^1$ | $\hat{\epsilon}_{3,t}^1$ | $\pi_{1,t}^{u,1}$ | $\pi_{3,t}^{u,1}$ | $\hat{\epsilon}_{1,t}^2$ | $\hat{\epsilon}_{3,t}^2$ | $\pi_{1,t}^{u,2}$ | $\pi_{3,t}^{u,2}$ |
| 21 | 31.15 | 8.31 | 14.87 | 14.87 | −31.15 | 8.31 | −17.67 | 1.77 |

TABLE VII
LMP AND UMP AT HOUR 21 ($\Lambda = 1, \Lambda^\Delta = 2$)

| Price | Bus1 | Bus2 | Bus3 | Bus4 | Bus5 | Bus6 |
|---|---|---|---|---|---|---|
| $\pi^e$ | 14.97 | 32.64 | 34.4 | 43.71 | 41.94 | 35.26 |
| $\pi^{u,up}$ | 14.87 | 14.87 | 16.63 | 25.94 | 24.17 | 17.49 |
| $\pi^{u,down}$ | -17.67 | 0 | 0 | 0 | 0 | 0 |

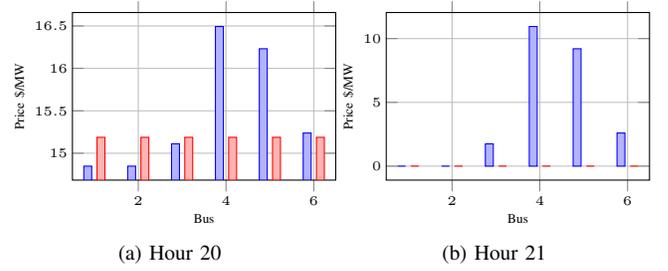

(a) Hour 20      (b) Hour 21

Fig. 3. Upward UMP (blue bar) and reserve price (red bar) with network constraint.

The load at Bus 4 has to pay the highest LMP $43.71/MWh. The UMPs are also different at various locations. The highest upward UMP at Hour 21 is also located at Bus 4. With these prices, the market participants can be paid and credited.

The LMP paid to G3 is $35.26/MWh on Bus 6, which is $17.49/MWh larger than its marginal cost. In the same time, the upward UMP is $17.49/MWh on Bus 6, which is exactly the difference between the LMP and G3's marginal cost. Hence, G3 is the UMP setter on Bus 6. The UMPs provide important price signals on the planning of renewable energy sources and storages. For example, the UMP at Bus 2 is relatively small, so it is an ideal location for renewable energy sources in terms of payment for uncertainties. In contrast, the UMP at Bus 4 is large, which may attract the long-term investment for storages or generation plants with large ramping rates.

*2) Comparison with Existing LMPs and Reserve Prices:*
The motivation of this part is to compare the proposed clearing scheme with the existing one. However, as the reserve is not robust in the traditional scheme, we cannot compare them fairly. With the observation that the transmission constraints are the most challenging one in the robust UC framework, we drop these constraints in this subsection and add reserve constraints as follows.

$$I_{i,t} P_i^{min} \leq Q_{i,t}^{down} + P_{i,t}, \quad Q_{i,t}^{up} + P_{i,t} \leq I_{i,t} P_i^{max}, \forall i, t \quad (24a)$$

$$-R_i^d I_{i,t} \Delta T \leq Q_{i,t}^{down}, \quad Q_{i,t}^{up} \leq R_i^u I_{i,t} \Delta T, \forall i, t \quad (24b)$$

$$\sum_i Q_{i,t}^{down} \leq \underline{R}_t, \quad \sum_i Q_{i,t}^{up} \leq \bar{R}_t, \forall t, \quad (24c)$$

where $Q_{i,t}^{up}$ and $Q_{i,t}^{down}$ are the largest upward and downward reserves, respectively. $\underline{R}_t$ and $\bar{R}_t$ are system-wide reserve requirements. Refer to [2], [30] for more details on the reserve formulations. In the experiment, $\Delta T$ is set to 1 and $\Lambda = 0.8, \Lambda^\Delta = 2$. The reserve requirements $\underline{R}_t$ and $\bar{R}_t$ are set to the lower and upper bounds of the system-wide uncertainty in (23b), respectively.

The results are as expected. The optimal solutions of the RUC and the standard UC with explicit reserve constraints are the same. LMPs calculated in the RUC and UC also have the same values. The UMPs calculated in the proposed mechanism are also exactly the same as the reserve prices in standard UC. Two things are verified with these results. First, without transmission constraints, the solution to standard UC can easily be robust by adding reserve constraints. Second, the proposed LMPs and UMPs are consistent with LMPs and reserve prices in the existing market when the transmission constraints are dropped.

When considering transmission constraints, the generation reserve cannot be guaranteed at bus levels in the traditional UC model. For simplicity, we assume that the 6 buses are in a zone. Consider the case where $\Lambda = 0.8, \Lambda^\Delta = 2$. The upward UMP and reserve price at Hour 20 are depicted in Fig. 3a. It is observed UMPs at Bus 1 and Bus 2 are lower than the traditional reserve prices. In the same time, the UMPs at Bus 4 and Bus 5 is higher than the traditional reserve prices. The differences are caused by the congestion of Line 1-4 for reserve delivery. It is worth mentioning that the LMP differences in two models are within 1% at Hour 20. The prices illustrated in Fig. 3b reveals another trends that the UMP may be higher than the traditional reserve prices. At Hour 21, the zonal reserve price is 0 while the UMPs are non-zeros at bus 4, 5, and 6. Because the constraint related with reserves in the RUC is stronger than the one in traditional UC model. Consequently, more expensive resources are used in RUC, which also generally leads to higher UMPs.

*3) FTR Underfunding:* When $\Lambda^\Delta = 2, \Lambda = 1$, the generation schedules at Hour 21 are 195.193MW, 25.577MW, and 16.54MW. The power flow of Line 2 is 97.63MW, which is 2.47MW smaller than its physical limit of 100MW. The transmission reserve 2.47 MW is kept to guarantee the delivery of the generation reserve. The binding constraint for Line 2 causes LMP differences. Hence, the FTR holder gets credits. Consider a set of FTR amounts [202.3429, 23.2771, −55.772, −94.924, −94.924, 20]. It can be verified that the FTR amounts satisfy the SFT in the FTR market. Then the total credit for the FTR holders is $5,554.77. However, the congestion cost in the DAM is $5,422.87. It means that the LMP congestion cost collected is not enough to cover the FTR credit. The FTR underfunding value is $5554.77 − $5422.87 = $131.90.

The revenue residues after UMP settlement is $131.9. It exactly covers the FTR underfunding in this scenario. Therefore, the revenue is adequate at Hour 21.

*B. IEEE 118-Bus System*

The simulations are performed for the IEEE 118-bus system with 54 thermal units and 186 branches in this section. The peak load is 6600MW. The detailed data including generator parameters, line reactance and ratings, and load profiles can be found at http://motor.ece.iit.edu/Data/RUC118UMP.xls. Two cases are studied in this section.

1) The uncertainty levels and load levels are changed to analyze the simulation results in the system level. The impact of transmission line capacity on prices is also studied.

TABLE VIII
OPERATION COST AND UMP PAYMENT ($, $\Lambda^\Delta = 10$)

| $\Lambda$ | Op. Cost | Un. Payment | Gen. Res. Credit | Rev. Res. |
|-----|-----|-----|-----|-----|
| 0.2 | 1,866,023 | 11,043 | 10,560 | 483 |
| 0.25 | 1,871,364 | 20,044 | 19,209 | 835 |
| 0.3 | 1,877,471 | 30,879 | 28,658 | 2,221 |

2) An energy storage is installed at a specified bus with high UMP to show the potential application of UMPs.

*1) Case 1:* We assume that the uncertainty sources are located at buses $(11, 15, 49, 54, 56, 59, 60, 62, 80, 90)$. The budget parameter $\Lambda^\Delta$ is set to 10 in this section. The bus-level uncertainty budget parameter $\Lambda$ changes from 0.2 to 0.3, and the bound of the uncertainty is the base load. The simulation results are shown in Table VIII. It can be observed that the total operation cost increases with increasing $\Lambda$. It indicates that a larger uncertainty level may increase the operation cost. The columns "Un. Payment" and "Gen. Res. Credit" denote the total payment from uncertainty sources and credit to generation reserves, respectively. The lowest payment is \$11,043 and the highest one is \$30,879. On the other hand, the credit entitled to the generation reserves is also a monotonically increasing function of $\Lambda$. When $\Lambda = 0.3$, the generation reserves have the highest credit. The last column "Rev. Res." shows that the revenue residues related to UMPs. It can be observed that the residue is always positive.

Fig. 4a in the next page depicts the heat map for the upward UMPs from Bus 80 to Bus 100 in 24 hours. The x-axis represents time intervals and the y-axis represents bus numbers. The color bar on the right shows different colors for various UMP values. For example, the \$0/MWh is denoted by the blue color at the bottom, and the \$18/MWh is represented by the dark red color on the top of the color bar. It can be observed that the uncertainty sources have system-wide unique UMPs at some intervals, such as Hours 8, 13, 15, and so on. It indicates that there is no transmission reserve in these hours. On the other hand, the UMPs at Hour 11 vary dramatically with different locations. The highest upward UMP is around \$18/MWh, and the lowest one is around \$2/MWh. According to the data shown in Fig. 4a, the high UMP at Bus 94 may attract investment of flexible resources, such as energy storages, in terms of generation reserve credit, and Bus 100 is an attractive location for the investment of renewable energy sources in terms of uncertainty payments.

Fig. 5 shows the uncertainty payment and generation reserve credit with respect to load levels. The base load level is set at 100%. Higher loads in general lead to more uncertainty payments and generation reserve credits. It is also consistent with the heat map of UMPs in Fig. 4a, where UMPs at peak load hours are high. It suggests that the generation reserves also become scarce resources when load levels are high.

The transmission line capacity plays an important role in the price calculation. Fig. 6 shows the LMPs and upward UMPs at Hour 11 with respect to increasing capacity of Line 94-100. The prices at Buses 88, 94, and 100 are depicted. When the line capacity increases from 165MW to 175MW, LMP at

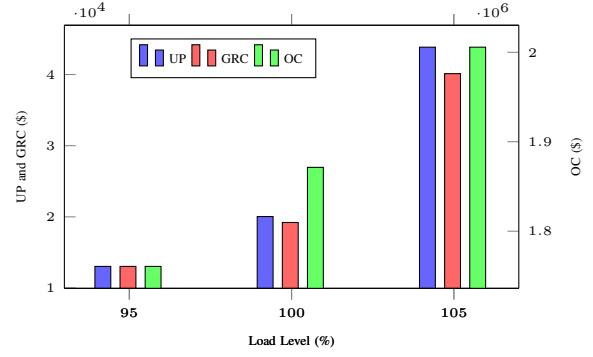

Fig. 5. Uncertainty payment (UP), generation reserve Credit (GRC), and operation cost (OC) with different load levels ($\Lambda = 0.25$)

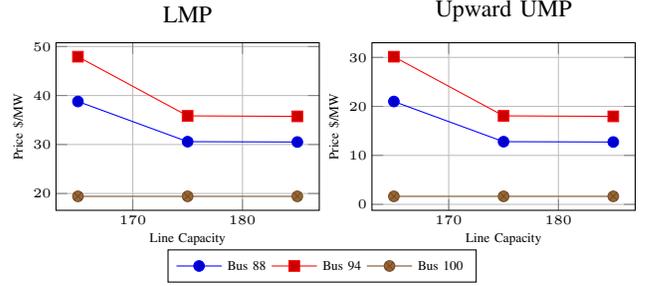

Fig. 6. LMP (left) and upward UMP (right) at Hour 11 with respect to increasing capacity of Line 94-100

Bus 94 decreases from \$47.92/MWh to \$35.84/MWh and that at Bus 88 also drops to \$30.58/MWh from \$38.78/MWh. The upward UMPs at Bus 94 and Bus 88 also drop by \$8.20/MWh and \$12.08/MWh, respectively. In contrast, the LMP and upward UMP at Bus 100, which is connected to Line 94-100, remain at \$19.42/MWh and \$1.64/MWh, respectively. It shows that the change of line capacity may only have impacts on the prices at some buses. When the line capacity further increases to 185MW from 175MW, the changes of LMPs and UMPs at Bus 94 and Bus 88 are within \$0.1/MWh, and there is still no change at Bus 100. It means that the additional 10MW cannot help deliver cheaper energy and reserves to Bus 94 and Bus 88. These results are also consistent with the analysis of the traditional LMPs [31].

*2) Case 2:* As discussed in Case 1, the upward UMP on Bus 94 is high at Hour 11. Assume that an energy storage (8MW/30MWh) is installed at Bus 94. A simple model for the energy storage is formulated as follows.

$$E_t = E_{t-1} + \rho^d P_t^D + \rho^c P_t^C, \forall t$$
$$0 \le E_t \le E^{\max}, \forall t$$
$$0 \le -P_t^D \le I_t^D R^D, \forall t$$
$$0 \le P_t^C \le I_t^C R^C, \forall t$$
$$I_t^D + I_t^C \le 1, \forall t$$
$$E_{N_T} = E_0,$$

where $E_t$ denotes the energy level, $P_t^D$ and $P_t^C$ represent the discharging and charging rates, and $I_t^D$ and $I_t^C$ are the indicators of discharging and charging. As the UMP is the major concern in this section, we use simplified parameters for storage. The discharging efficiency $\rho^d$ and charging efficiency $\rho^c$ are set to 100%. The capacity $E^{\max}$ and initial energy level $E_0$ are set to 30 MWh and 15 MWh, respectively. The maximal charging rate $R^D$ and discharging rate $R^C$ are set to 8 MW/h.

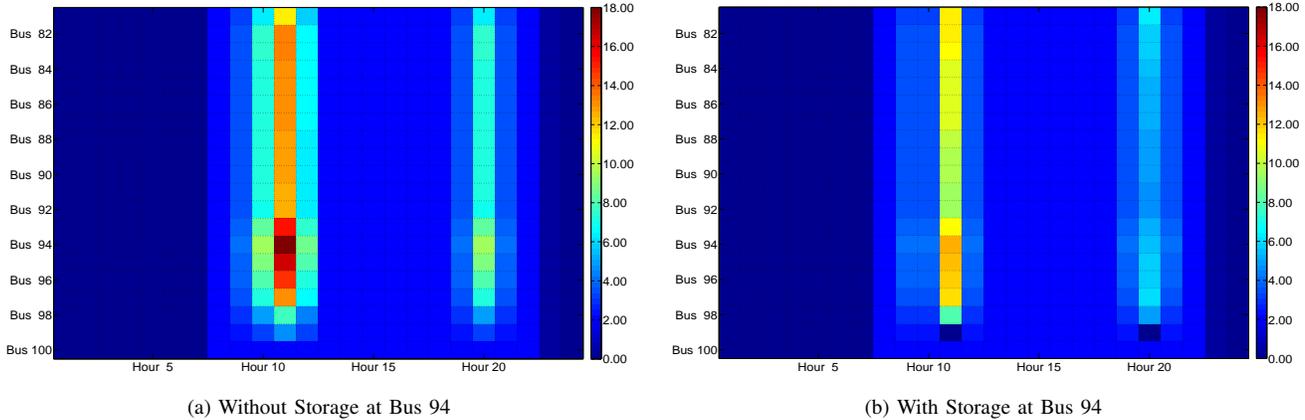

(a) Without Storage at Bus 94

(b) With Storage at Bus 94

Fig. 4. Heat Map for Upward UMPs ($\Lambda = 0.3$). Different colors represent various UMPs. Figure (a) depicts the UMPs from Bus 80 to Bus 100 in 24 hours without the energy storage at Bus 94. Figure (b) depicts the new UMPs after the energy storage is sited at Bus 94.

By siting the energy storage, we can lower the new operation cost to \$1,875,211 from \$1,877,471. The payment collected from the uncertainty sources becomes \$27,473, and the credit to generation reserves decreases to \$24,289. Compared to the data in Table VIII, the energy storage also helps to reduce the payment related to UMPs. The storage is entitled to \$1326 generation reserve credit. Fig. 4b depicts the new upward UMPs after the installation of the energy storage. Compared to that in Fig. 4a, the upward UMP for Hour 11 at Bus 94 decreases a lot. The UMPs for Hour 10 and 12 are also lower. It suggests that siting the energy storage at Bus 94 effectively lower the generation reserve price.

The simulation results demonstrate that flexible resources can lower the UMPs, and UMPs provide the investment signal at locations where generation reserves are scarce resources.

## IV. Conclusions

A novel market model in this paper clears uncertainty, energy, and generation reserve simultaneously within the RUC framework in DAM. The uncertainty sources are charged and the generator reserve providers are credited based on the proposed UMP. The UMP formulation is derived within a robust optimization framework. We also characterize the market equilibrium for the new market clearing mechanism. As the market clearing mechanism is established within the robust optimization framework, the robustness of the dispatch is guaranteed. The optimal reserves for uncertainty accommodation are obtained in the model. The UMP proposed in this paper can effectively address the issue on how to charge and credit the uncertainties and generation reserve fairly in the market with RES.

Our study also shows that traditional pricing mechanism within RUC framework may lead to FTR underfunding. The proposed market clearing mechanism can address this issue. Our study shows load serving entities can have lower energy prices within the new market scheme, as the reserve fees are paid by uncertainty sources.

Many potential applications on UMPs are open. As UMPs are unified prices of uncertainties and reserves, it is interesting to investigate the optimal strategy for the one who is the uncertainty source as well as reserve provider in the market (e.g. the wind generation company with energy storages). The UMPs derived in this paper also provide an important price signal for the long-term investment of flexible resources. When the upward UMP or downward UMP at a bus is high, the investor can get more return in terms of generation reserves.

Another potential future research on UMP is to study how to determine the budget uncertainty set in the market. Modeling the traditional spinning reserve for the contingency [13], [30] is also our future work. In this paper, the demand bids are not considered. We have forecasted load, forecasted RES, and uncertainty of load and RES for market clearing with a single-sided model. In an extended double-sided model, we can have demand bids, forecasted RES, and uncertainty of load and RES for market clearing. The forecasted load, forecasted RES, uncertainty of load and RES can be used in RAC.

## Appendix A
## Detailed Formulation for Problem (RUC)

$$\text{(RUC)} \min_{(x,y,z,I,P) \in \mathcal{F}} \sum_t \sum_i \left( C_i^P(P_{i,t}) + C_i^I(I_{i,t}) \right) \quad (25a)$$

$$\text{s.t.} \quad \sum_i P_{i,t} = \sum_m d_{m,t}, \forall t. \quad (25b)$$

$$-F_l \le \sum_m \Gamma_{l,m} \left( \sum_{i \in \mathcal{G}(m)} P_{i,t} - d_{m,t} \right) \le F_l, \forall l, t \quad (25c)$$

$$I_{i,t} P_i^{\min} \le P_{i,t} \le I_{i,t} P_i^{\max}, \forall i, t \quad (26a)$$

$$P_{i,t} - P_{i,(t-1)} \le r_i^u(1 - y_{i,t}) + P_i^{\min} y_{i,t}, \forall i, t \quad (26b)$$

$$-P_{i,t} + P_{i,(t-1)} \le r_i^d(1 - z_{i,t}) + P_i^{\min} z_{i,t}, \forall i, t \quad (26c)$$

minimum on/off time limit

and

$$\mathcal{F} := \Big\{ (x, y, z, I, P) : \forall \epsilon \in \mathcal{U}, \exists \Delta P \text{ such that}$$

$$\sum_i \Delta P_{i,t} = \sum_m \epsilon_{m,t}, \forall t, \tag{27a}$$

$$I_{i,t} P_i^{\min} \leq P_{i,t} + \Delta P_{i,t} \leq I_{i,t} P_i^{\max}, \forall i, t \tag{27b}$$

$$-R_i^d(1 - z_{i,t+1}) \leq \Delta P_{i,t} \leq R_i^u(1 - y_{i,t}), \forall i, t \tag{27c}$$

$$\Delta P_{m,t}^{\text{inj}} = \sum_{i \in \mathcal{G}(m)} \Delta P_{i,t} - \epsilon_{m,t}, \forall m, t \tag{27d}$$

$$-F_l \leq \sum_m \Gamma_{l,m}(P_{m,t}^{\text{inj}} + \Delta P_{m,j}^{\text{inj}}) \leq F_l, , \forall l, t \Big\}. \tag{27e}$$

The basic idea of the above model is to find a robust UC and dispatch for the base-case scenario. In the base-case scenario, (25b) denotes the load balance constraint; (25c) represents the transmission line constraint; (26a) denotes the unit capacity limit constraint; (26b)-(26c) denote the unit ramping up/down limits. $I_{i,t}$, $y_{i,t}$, and $z_{i,t}$ are the indicators of the unit being on, started-up, and shutdown, respectively. Units also respect the minimum on/off time constraints which are related to these binary variables [1]. The UC and dispatch solution are immunized against any uncertainty $\epsilon \in \mathcal{U}$. When uncertainty $\epsilon$ occurs, it is accommodated by the generation adjustment $\Delta P_{i,t}$ (27a). Generation dispatch is also enforced by the capacity limits (27b). (27c) models the ramping rate limits of generation adjustment $\Delta P_{i,t}$. In fact, the right and left hand sides of (27c) can correspond to a response time $\Delta T$, which is similar to the 10-min or 30-min reserves in the literatures [30]. (27e) stands for the network constraint after uncertainty accommodation.

## Appendix B
## Detailed Formulation for Problem (MP) and (SP)

$$\text{(MP)} \min_{(x,y,z,I,P,\Delta P)} \sum_t \sum_i \Big( C_i^P(P_{i,t}) + C_i^I(I_{i,t}) \Big)$$

S.T.  (25b), (25c), (26a)-(26c), minimum on/off time limit

$$\sum_i \Delta P_{i,t}^k = \sum_m \epsilon_{m,t}^k, \forall t, \forall k \in \mathcal{K} \tag{28a}$$

$$I_{i,t} P_i^{\min} \leq P_{i,t} + \Delta P_{i,t}^k \leq I_{i,t} P_i^{\max}, \forall i, t, \forall k \in \mathcal{K} \tag{28b}$$

$$\Delta P_{i,t}^k \leq R_i^u(1 - y_{i,t}), \forall i, t, \forall k \in \mathcal{K} \tag{28c}$$

$$-\Delta P_{i,t}^k \leq R_i^d(1 - z_{i,t+1}), \forall i, t, \forall k \in \mathcal{K} \tag{28d}$$

$$-F_l \leq \sum_m \Gamma_{l,m}(P_{m,t}^{\text{inj}} + \Delta P_{m,t}^{\text{inj},k}) \leq F_l, \forall k \in \mathcal{K}, \forall l, t \tag{28e}$$

$$\Delta P_{m,t}^{\text{inj},k} = \sum_{i \in \mathcal{G}(m)} \Delta P_{i,t}^k - \epsilon_{m,t}^k, \forall m, t, \forall k \in \mathcal{K}, \tag{28f}$$

and

$$\text{(SP)} \max_{\epsilon \in \mathcal{U}} \min_{(s^+, s^-, \Delta P) \in \mathcal{R}(\epsilon)} \sum_m \sum_t (s_{m,t}^+ + s_{m,t}^-) \tag{29a}$$

$$\mathcal{R}(\epsilon) := \Big\{ (s^+, s^-, \Delta P) : \tag{29b}$$

$$\sum_i \Delta P_{i,t} = \sum_m (\epsilon_{m,t} + s_{m,t}^+ - s_{m,t}^-), \forall m, t \tag{29c}$$

$$-F_l \leq \sum_m \Gamma_{l,m}\Big(P_{m,t}^{\text{inj}} + \Delta P_{m,t}^{\text{inj}}\Big) \leq F_l, \forall l, t \tag{29d}$$

$$\Delta P_{m,t}^{\text{inj}} = \sum_{i \in \mathcal{G}(m)} \Delta P_{i,t} - (\epsilon_{m,t} + s_{m,t}^+ - s_{m,t}^-) \tag{29e}$$

$$s_{m,t}^+, s_{m,t}^- \geq 0, \forall m, t \tag{29f}$$

$$(27b), (27c) \Big\}$$

where $\mathcal{K}$ is the index set for uncertainty points $\hat{\epsilon}$ which are dynamically generated in (SP) with iterations. It should be noted that $\hat{\epsilon}^{\boldsymbol{k}}$ is the extreme point of $\mathcal{U}$. Variable $\Delta P_{i,t}^k$ is associated with $\hat{\epsilon}^{\boldsymbol{k}}$. The objective function in (SP) is the summation of non-negative slack variables $s_{m,t}^+$ and $s_{m,t}^-$, which evaluates the violation associated with the solution $(x, y, z, I, P)$ from (MP). $s_{m,t}^+$ and $s_{m,t}^-$ are also explained as un-followed uncertainties (generation or load shedding) due to system limitations.

## Appendix C
## Lagrangian Function for Problem (RED)

Please check equation (30) in the next page.

## Appendix D
## Proofs for Lemmas and Theorems

### A. Proof of Lemma 1

*Proof.* Consider $\hat{\epsilon}_{m,t}^k > 0, \pi_{m,t}^{u,k} < 0$. With a small perturbation $\delta > 0$ to $\hat{\epsilon}_{m,t}^k$, we replace $\hat{\epsilon}_{m,t}^k$ with $\hat{\epsilon}_{m,t}^k - \delta$ in (RED). As the $\pi_{m,t}^{u,k} < 0$, then the optimal value to problem (RED) increases. It means that there are violations for the original optimal solution $P_{i,t}$ to problem (RED) with $\hat{\epsilon}_{m,t}^k - \delta$. Hence, the optimal solution $P_{i,t}$ to problem (RED) cannot be immunized against the uncertainty $\hat{\epsilon}_{m,t}^k - \delta$. It contradicts with the robustness of the solution $P_{i,t}$. Therefore, if $\hat{\epsilon}_{m,t}^k > 0$, then $\pi_{m,t}^{u,k} \geq 0$. Similarly, if $\hat{\epsilon}_{m,t}^k < 0$, then $\pi_{m,t}^{u,k} \leq 0$. $\qquad \square$

### B. Proof of Lemma 2

*Proof.* Assume $i \in \mathcal{G}(m)$, according to the KKT condition

$$\frac{\partial \mathcal{L}(P, \Delta P, \lambda, \alpha, \beta, \eta)}{\partial \Delta P_{i,t}^k} = 0 \tag{31}$$

at the optimal point, we have

$$\bar{\beta}_{i,t}^k - \underline{\beta}_{i,t}^k + \bar{\alpha}_{i,t}^k - \underline{\alpha}_{i,t}^k - \lambda_t^k + \sum_l (\bar{\eta}_{l,t}^k - \underline{\eta}_{l,t}^k) \Gamma_{l,m} = 0. \tag{32}$$

Then (17) holds. If $\pi_{m,t}^{u,k} > 0$, then $\bar{\beta}_{i,t}^k + \bar{\alpha}_{i,t}^k > 0$ as $\bar{\beta}_{i,t}^k, \underline{\beta}_{i,t}^k, \bar{\alpha}_{i,t}^k$, and $\underline{\alpha}_{i,t}^k$ are non-negative. According to the complementary conditions for (7b) and (7d), at least one of (7b) and (7d) is binding. Hence,

$$\Delta P_{i,t}^k = \min \Big\{ \hat{I}_{i,t} P_i^{\max} - P_{i,t}, R_i^u(1 - \hat{y}_{i,t}) \Big\}$$

holds. Similarly, the other equation holds when $\pi_{m,t}^{u,k} < 0$. $\qquad \square$

$$
\begin{aligned}
& \mathcal{L}(P, \Delta P, \lambda, \alpha, \beta, \eta) \\
= & \sum_t \sum_i C_i^P(P_{i,t}) + \sum_t \lambda_t \Big( \sum_m d_{m,t} - \sum_i P_{i,t} \Big) + \sum_t \sum_i \Big( \bar{\beta}_{i,t}(P_{i,t} - \hat{I}_{i,t} P_i^{\max}) + \underline{\beta}_{i,t}(\hat{I}_{i,t} P_i^{\min} - P_{i,t}) \Big) \\
& + \sum_t \sum_i \Big( \bar{\alpha}_{i,t} \Big( P_{i,t} - P_{i,t-1} - r_i^u(1 - \hat{y}_{i,t}) - P_{i,t}^{\min} \hat{y}_{i,t} \Big) + \underline{\alpha}_{i,t} \Big( P_{i,t-1} - P_{i,t} - r_i^d(1 - \hat{z}_{i,t}) - P_{i,t}^{\min} \hat{y}_{i,t} \Big) \Big) \\
& + \sum_t \sum_l \Big( \bar{\eta}_{l,t} \Big( \sum_m \Gamma_{l,m} P_{m,t}^{\mathrm{inj}} - F_l \Big) - \underline{\eta}_{l,t} \Big( \sum_m \Gamma_{l,m} P_{m,t}^{\mathrm{inj}} + F_l \Big) \Big) + \sum_{k \in \mathcal{K}} \sum_t \lambda_t^k \Big( \sum_m \epsilon_{m,t}^k - \sum_i \Delta P_{i,t}^k \Big) \\
& + \sum_{k \in \mathcal{K}} \sum_t \sum_i \Big( \bar{\beta}_{i,t}^k (P_{i,t} + \Delta P_{i,t}^k - \hat{I}_{i,t} P_i^{\max}) + \underline{\beta}_{i,t}^k (\hat{I}_{i,t} P_i^{\min} - P_{i,t} - \Delta P_{i,t}^k) \Big) \\
& + \sum_{k \in \mathcal{K}} \sum_t \sum_i \Big( \bar{\alpha}_{i,t}^k \Big( \Delta P_{i,t}^k - R_i^u(1 - \hat{y}_{i,t}) \Big) - \underline{\alpha}_{i,t}^k \Big( \Delta P_{i,t}^k + R_i^d(1 - \hat{z}_{i,t}) \Big) \Big) \\
& + \sum_{k \in \mathcal{K}} \sum_t \sum_l \Big( \bar{\eta}_{l,t}^k \Big( \sum_m \Gamma_{l,m}(P_{m,t}^{\mathrm{inj}} + \Delta P_{m,t}^{\mathrm{inj},k}) - F_l \Big) - \underline{\eta}_{l,t}^k \Big( \sum_m \Gamma_{l,m}(P_{m,t}^{\mathrm{inj}} + \Delta P_{m,t}^{\mathrm{inj},k}) + F_l \Big) \Big)
\end{aligned}
\tag{30}
$$

---

### C. Proof of Theorem 1

*Proof.* The energy congestion cost at $t$ is

$$
\sum_m \Big( \pi_{m,t}^{\mathrm{e}} d_{m,t} - \sum_{i \in \mathcal{G}(m)} \pi_{m,t}^{\mathrm{e}} P_{i,t} \Big) \tag{33}
$$
$$
= \sum_m \pi_{m,t}^{\mathrm{e}} P_{m,t}^{\mathrm{inj}}
$$
$$
= \sum_m \sum_l \Gamma_{l,m} \Big( \bar{\eta}_{l,t} + \sum_{k \in \mathcal{K}} \bar{\eta}_{l,t}^k - \underline{\eta}_{l,t} - \sum_{k \in \mathcal{K}} \underline{\eta}_{l,t}^k \Big) P_{m,t}^{\mathrm{inj}}
$$
$$
= \sum_l \Big( \bar{\eta}_{l,t} + \sum_{k \in \mathcal{K}} \bar{\eta}_{l,t}^k \Big) \Big( F_l - \Delta f_{l,t}^{\mathrm{pos}} \Big)
$$
$$
- \sum_l \Big( \underline{\eta}_{l,t} + \sum_{k \in \mathcal{K}} \underline{\eta}_{l,t}^k \Big) \Big( \Delta f_{l,t}^{\mathrm{neg}} - F_l \Big)
$$
$$
= \sum_l \Big( \bar{\eta}_{l,t} + \sum_{k \in \mathcal{K}} \bar{\eta}_{l,t}^k + \underline{\eta}_{l,t} + \sum_{k \in \mathcal{K}} \underline{\eta}_{l,t}^k \Big) F_l
$$
$$
- \sum_l \sum_{k \in \mathcal{K}} \Big( \bar{\eta}_{l,t}^k \Delta f_{l,t}^{\mathrm{pos}} + \underline{\eta}_{l,t}^k \Delta f_{l,t}^{\mathrm{neg}} \Big)
$$
$$
- \sum_l \Big( \bar{\eta}_{l,t} \Delta f_{l,t}^{\mathrm{pos}} + \underline{\eta}_{l,t} \Delta f_{l,t}^{\mathrm{neg}} \Big)
$$
$$
= \sum_l \Big( \bar{\eta}_{l,t} + \sum_{k \in \mathcal{K}} \bar{\eta}_{l,t}^k + \underline{\eta}_{l,t} + \sum_{k \in \mathcal{K}} \underline{\eta}_{l,t}^k \Big) F_l
$$
$$
- \sum_l \sum_{k \in \mathcal{K}} \Big( \bar{\eta}_{l,t}^k \Delta f_{l,t}^{\mathrm{pos}} + \underline{\eta}_{l,t}^k \Delta f_{l,t}^{\mathrm{neg}} \Big)
$$

The first equality holds following the definition of net power injection. The second equality holds according to (8) and $\sum_m P_{m,t}^{\mathrm{inj}} = 0$. The third equality holds following (19) and (20). The sign change of $\underline{\eta}_{l,t}$ and $\sum_l \underline{\eta}_{l,t}^k$ in the third equality is because of the definition of power flow direction. According to the complementary conditions, the third term in the fourth equality must be zero based on the following three cases.

1) If $\bar{\eta}_{l,t} \neq 0$, then $\Delta f_{l,t}^{\mathrm{pos}} = 0$, and $\underline{\eta}_{l,t} = 0$.
2) If $\bar{\eta}_{l,t} = 0$ and $\underline{\eta}_{l,t} \neq 0$, then $\Delta f_{l,t}^{\mathrm{neg}} = 0$.

3) $\bar{\eta}_{l,t} = 0$ and $\underline{\eta}_{l,t} = 0$.

The second term in the last equality corresponds to (21).

The credits to FTR holders can be written as

$$
\sum_{(m \to n)} (\pi_{m,t}^{\mathrm{e}} - \pi_{n,t}^{\mathrm{e}}) \mathrm{FTR}_{m \to n} \tag{34}
$$
$$
= \sum_{(m \to n)} \begin{pmatrix} \lambda_t - \sum_l \Gamma_{l,m}(\bar{\eta}_{l,t} - \underline{\eta}_{l,t}) \\ - \sum_l \sum_{k \in \mathcal{K}} \Gamma_{l,m}(\bar{\eta}_{l,t}^k - \underline{\eta}_{l,t}^k) \\ -\lambda_t + \sum_l \Gamma_{l,n}(\bar{\eta}_{l,t} - \underline{\eta}_{l,t}) \\ + \sum_l \sum_{k \in \mathcal{K}} \Gamma_{l,n}(\bar{\eta}_{l,t}^k - \underline{\eta}_{l,t}^k) \end{pmatrix} \mathrm{FTR}_{m \to n}
$$
$$
= \sum_{(m \to n)} \sum_l \begin{pmatrix} (\Gamma_{l,n} - \Gamma_{l,m})(\bar{\eta}_{l,t} + \sum_{k \in \mathcal{K}} \bar{\eta}_{l,t}^k) \\ -(\Gamma_{l,n} - \Gamma_{l,m})(\underline{\eta}_{l,t} + \sum_{k \in \mathcal{K}} \underline{\eta}_{l,t}^k) \end{pmatrix} \mathrm{FTR}_{m \to n}
$$
$$
\leq \sum_l \Big( \bar{\eta}_{l,t} + \sum_{k \in \mathcal{K}} \bar{\eta}_{l,t}^k + \underline{\eta}_{l,t} + \sum_{k \in \mathcal{K}} \underline{\eta}_{l,t}^k \Big) F_l
$$

The first equality holds according to (8). The inequality is true as the amount of $\mathrm{FTR}_{m \to n}$ respects

$$
-F_l \leq \sum_{m \to n} (\Gamma_{l,m} - \Gamma_{l,n}) \mathrm{FTR}_{m \to n} \leq F_l.
$$

according to the SFT for FTR market [18], [20], [21]. The right-hand-side of the inequality is the first term in the last equality of (33). Based on (33) and (34), the maximum difference between the FTR credit and the energy congestion cost is equal to the transmission reserve credit. That is, the maximum FTR underfunding is (21). $\qquad \square$

### D. Proof of Theorem 2

*Proof.* According to Theorem 1, the FTR underfunding value is (21) due to the deficiency of the energy congestion cost.

Therefore, we need to prove that the money collected from uncertainty sources can cover the FTR underfunding and credits to generation reserve.

Without loss of generality, we consider the payment collected from uncertainty sources at time $t$ for $\hat{\epsilon}^k$

$$
\begin{aligned}
&\sum_m \pi_{m,t}^{u,k} \hat{\epsilon}_{m,t}^k \\
&= \sum_m \left( \lambda_t^k - \sum_l \Gamma_{l,m} \left( \bar{\eta}_{l,t}^k - \underline{\eta}_{l,t}^k \right) \right) \hat{\epsilon}_{m,t}^k \\
&= \sum_i \Delta P_{i,t}^k \lambda_t^k - \sum_m \sum_l \Gamma_{l,m} \left( \bar{\eta}_{l,t}^k - \underline{\eta}_{l,t}^k \right) \left( \sum_{i \in \mathcal{G}(k)} \Delta P_{i,t}^k \right) \\
&\quad + \sum_l \left( \bar{\eta}_{l,t}^k \Delta f_{l,t}^{\text{pos}} + \underline{\eta}_{l,t}^k \Delta f_{l,t}^{\text{neg}} \right) \\
&= \sum_m \sum_{i \in \mathcal{G}(m)} \pi_{m,t}^{u,k} \Delta P_{i,t}^k + \sum_l \left( \bar{\eta}_{l,t}^k \Delta f_{l,t}^{\text{pos}} + \underline{\eta}_{l,t}^k \Delta f_{l,t}^{\text{neg}} \right) \\
&= \sum_i \pi_{m_i,t}^{u,k} \Delta P_{i,t}^k + \sum_l \left( \bar{\eta}_{l,t}^k \Delta f_{l,t}^{\text{pos}} + \underline{\eta}_{l,t}^k \Delta f_{l,t}^{\text{neg}} \right)
\end{aligned}
$$

The first equality holds according to (9). According to (7a), (7f), and (7g), the $\sum_m \sum_l \Gamma_{l,m} \bar{\eta}_{l,t}^k \hat{\epsilon}_{m,t}^k$ in the second line can be rewritten as

$$
\begin{aligned}
&\sum_m \sum_l \Gamma_{l,m} \bar{\eta}_{l,t}^k \left( \sum_{i \in \mathcal{G}(m)} (\Delta P_{i,t}^k + P_{i,t}) - d_{m,t} \right) - \sum_l \bar{\eta}_{l,t}^k F_l \\
&= \sum_m \sum_l \Gamma_{l,m} \bar{\eta}_{l,t}^k \sum_{i \in \mathcal{G}(m)} \Delta P_{i,t}^k + \sum_l \bar{\eta}_{l,t}^k \left( \sum_m \Gamma_{l,m} P_{m,t}^{\text{inj}} - F_l \right) \\
&= \sum_m \sum_l \Gamma_{l,m} \bar{\eta}_{l,t}^k \left( \sum_{i \in \mathcal{G}(m)} \Delta P_{i,t}^k \right) - \sum_l \bar{\eta}_{l,t}^k \Delta f_{l,t}^{\text{pos}}
\end{aligned}
$$

$\sum_m \sum_l \Gamma_{l,m} \underline{\eta}_{l,t}^k \hat{\epsilon}_{m,t}^k$ can be reformulated similarly. Hence, the second equality holds. The third equality holds from (9). Therefore,

$$
\sum_m \sum_t \Psi_{m,t} = \sum_i \sum_t \Theta_{i,t}^G + \sum_l \sum_t \Theta_{l,t}^T
$$

holds. That is, the uncertainty payment covers the generation reserve credit and transmission reserve credit. Then, following the energy congestion cost shown in (33),

$$
\begin{aligned}
&\sum_m \sum_t \pi_{m,t}^e d_{m,t} + \sum_m \sum_t \Psi_{m,t} \\
&\geq \sum_i \sum_t \pi_{m_i,t}^e P_{i,t} + \sum_i \sum_t \Theta_{i,t}^G \\
&\quad + \sum_t \sum_{m \to n} (\pi_{m,t}^e - \pi_{n,t}^e) \text{FTR}_{m \to n}
\end{aligned}
$$

holds. That is, the total payments collected from loads and uncertainty sources can cover the total credits to energy, generation reserve, and FTR holders. So, the revenue adequacy of the proposed market clearing mechanism is guaranteed. □

### E. Proof of Competitive Equilibrium

*Proof.* $P_{i,t}$ and $(Q_{i,t}^{\text{up}}, Q_{i,t}^{\text{down}})$ are coupled by constraints (15) and (16). According to (17), we can rewrite generation reserve credit as

$$
\pi_{m,t}^{u,\text{up}} Q_{i,t}^{\text{up}} + \pi_{m,t}^{u,\text{down}} Q_{i,t}^{\text{down}} = \sum_{k \in \mathcal{K}} \pi_{m,t}^{u,k} \Delta P_{i,t}^k
$$

$$
\begin{aligned}
&= \sum_{k \in \mathcal{K}} (\bar{\beta}_{i,t}^k - \underline{\beta}_{i,t}^k + \bar{\alpha}_{i,t}^k - \underline{\alpha}_{i,t}^k) \Delta P_{i,t}^k \\
&= \sum_{k \in \mathcal{K}} \left( \begin{aligned} &\bar{\beta}_{i,t}^k (\hat{I}_{i,t} P_i^{\max} - P_{i,t}) + \underline{\beta}_{i,t}^k (P_{i,t} - \hat{I}_{i,t} P_i^{\min}) \\ &+ \bar{\alpha}_{i,t}^k (R_i^u (1 - \hat{y}_{i,t})) + \underline{\alpha}_{i,t}^k R_i^d (1 - \hat{z}_{i,t+1}) \end{aligned} \right) \quad (35)
\end{aligned}
$$

Substituting (35) into problem (PMP$_i$), we can decouple $P_{i,t}$ and $(Q_{i,t}^{\text{up}}, Q_{i,t}^{\text{down}})$. In fact, we also get all terms related to $P_{i,t}$ in Lagrangian $\mathcal{L}(P, \lambda, \alpha, \beta, \eta)$ for problem (RED). Since the problem (RED) is a linear programming problem, the saddle point $\hat{P}_{i,t}$, which is the optimal solution to (RED), is also the optimal solution to (PMP$_i$). Consequently, unit $i$ is not inclined to deviate its output level as it can obtain the maximum profit by following the ISO's dispatch instruction $\hat{P}_{i,t}$. Therefore, dispatch $\hat{P}_{i,t}$ and price signal $(\pi_{m,t}^e, \pi_{m,t}^{u,k})$ constitute a competitive partial equilibrium [27]. □

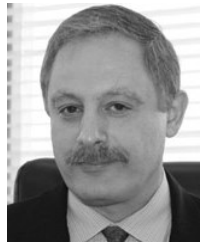

**Mohammad Shahidehpour** (F'01) received his Ph.D. degree from the University of Missouri in 1981 in electrical engineering. He is currently the Bodine Chair Professor and Director of the Robert W. Galvin Center for Electricity Innovation at the Illinois Institute of Technology, Chicago. He is the founding Editor-in-Chief of IEEE Transactions on Smart Grid. He is a member of US National Academy of Engineering (NAE).

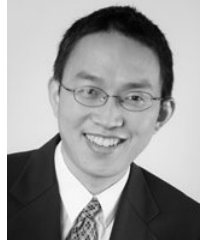

**Zuyi Li** (SM'09) received the B.S. degree from Shanghai Jiaotong University, Shanghai, China, in 1995, the M.S. degree from Tsinghua University, Beijing, China, in 1998, and the Ph.D. degree from the Illinois Institute of Technology (IIT), Chicago, in 2002, all in electrical engineering. Presently, he is a Professor in the Electrical and Computer Engineering Department at IIT. His research interests include economic and secure operation of electric power systems, cyber security in smart grid, renewable energy integration, electric demand management of data centers, and power system protection.

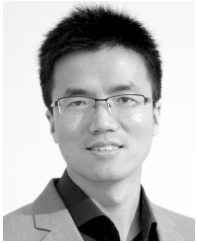

**Hongxing Ye** (S'14-m'16) received his B.S. degree in Information Engineering, in 2007, and M.S. degree in Systems Engineering, in 2011, both from Xi'an Jiaotong University, China, and the Ph.D. degree in Electrical Engineering from the Illinois Institute of Technology, Chicago in 2016. His research interests include large-scale optimization in power systems, electricity market, renewable integration, and cyber-physical system security in smart grid. He is "Outstanding Reviewer" for IEEE Transactions on Power Systems and IEEE Transactions on Sustainable Energy in 2015. He received Sigma Xi Research Excellence Award at Illinois Institute of Technology in 2016.

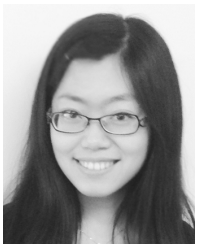

**Yinyin Ge** (S'14) received the B.S. degree (2008) in Automation and M.S. degree (2011) in Systems Engineering from Xian Jiaotong University, China. She also received Ph.D. degree (2016) in Electrical Engineering at Illinois Institute of Technology, Chicago. Her research interests are power system optimization and modeling; PMU applications in Smart Grid; monitoring, visualization, and state estimation for distribution systems.